\documentstyle[amscd,amssymb,verbatim,12pt]{amsart}
\newcommand{\gd}{\mbox{\rm geom dim} \:}
\newcommand{\Rad}{\mbox{\rm Rad}}
\newcommand{\SL}{\mbox{\rm SL}}
\newcommand{\GL}{\mbox{\rm GL}}
\newcommand{\Sl}{\mbox{\rm sl}}
\newcommand{\SP}{\mbox{\rm Sp}}
\newcommand{\SO}{\mbox{\rm SO}}
\newcommand{\SU}{\mbox{\rm SU}}
\newcommand{\PSL}{\mbox{\rm PSL}}
\newcommand{\Def}{\mbox{\rm def}}

\newcommand{\Ker}{\mbox{\rm Ker}}

\newcommand{\Aut}{\mbox{\rm Aut}}
\newcommand{\tr}{\mbox{\rm tr}}

\newcommand{\rank}{\mbox{\rm rank}}

\newcommand{\HH}{\mbox{\rm H}}

\newcommand{\vol}{\mbox{\rm vol}}

\newcommand{\R}{{\Bbb R}}
\newcommand{\C}{{\Bbb C}}

\newcommand{\Z}{{\Bbb Z}}

\theoremstyle{plain}

\newtheorem{lemma}{Lemma}
\newtheorem{theorem}{Theorem}

\numberwithin{equation}{section}
\renewcommand{\rm}{\normalshape}
\begin{document}
\title{Deficiencies of Lattice Subgroups of Lie Groups}
\author{John Lott}
\address{Department of Mathematics\\
University of Michigan\\
Ann Arbor, MI  48109-1109\\
USA}
\email{lott@@math.lsa.umich.edu}
\thanks{Research supported by NSF grant DMS-9704633}
\date{October 7, 1997}
\maketitle
\begin{abstract}
Let $\Gamma$ be a lattice in a connected Lie group.  We show that
besides a few exceptional cases, the deficiency of $\Gamma$ is nonpositive.
\end{abstract}
\section{Introduction}
If $\Gamma$ is a finitely presented group, its deficiency $\Def(\Gamma)$
is the maximum, over all finite presentations of $\Gamma$, of the
number of generators minus the number of relations.
If $G$ is a connected Lie group, a lattice in $G$ is a 
discrete subgroup  $\Gamma$
such that $G/\Gamma$ has finite volume. It is uniform if 
$G/\Gamma$ is compact.  Lubotzky proved the 
following result \cite[Proposition 6.2]{Lubotzky (1983)}: 
\begin{theorem} (Lubotzky)
Let $\Gamma$ be a lattice in a simple Lie group $G$. \\
(a) If $\R-\rank(G) \ge 2$ or $G = \SP(n,1)$ or $G= F_4$, then
$\Def(\Gamma) \le 0$. \\
(b) If $G = \SO(n,1)$ (for $n \ge 3$) or $G = \SU(n,1)$
(for $n \ge 2$), then
$\Def(\Gamma) \le 1$.
\end{theorem}

We give an improvement of Lubotzky's result.
\begin{theorem} \label{mainprop}
Let $G$ be a connected Lie group. Let $\Gamma$ be a lattice in
$G$. If $\Def(\Gamma) > 0$ then
1. $\Gamma$ has a finite normal subgroup $F$ such that $\Gamma/F$ is a 
lattice in $\PSL_2( \R)$, or\\
2. $\Def(\Gamma) = 1$ and either \\
A. $\Gamma$ is isomorphic to 
a torsion-free nonuniform lattice in $\R \times \PSL_2( \R)$ or
$\PSL_2( \C)$, or \\
B. $\Gamma$ is $\Z$, $\Z^2$ or the fundamental group of a Klein bottle.
\end{theorem}

The examples in case 2 do have deficiency one
\cite{Epstein (1961)}. A free group on $r$ generators, $r > 1$, has deficiency
$r$ and gives an example of case 1.

In some cases, we have sharper bounds on $\Def(\Gamma)$.
\begin{theorem} \label{improve}
1. If $\Gamma$ is a lattice in $\SO(4,1)$ then
\begin{equation}
\Def(\Gamma) \le 1 - \frac{3}{4 \pi^2} \: \vol(H^4/\Gamma).
\end{equation}
2. If $\Gamma$ is a lattice in $\SU(2,1)$ then
\begin{equation}
\Def(\Gamma) \le 1 - \frac{6}{\pi^2} \: \vol(\C H^2/\Gamma).
\end{equation}
(We normalize $\C H^2$ to have sectional curvatures between $-4$ and $-1$.)\\
3. If $\Gamma$ is a lattice in $\PSL_2(\R) \times \PSL_2( \R)$ then
\begin{equation}
\Def(\Gamma) \le 1 - \frac{1}{4 \pi^2} \: \vol((H^2 \times H^2)/\Gamma).
\end{equation}
\end{theorem}
\section{Proofs}
To prove Theorems \ref{mainprop} and \ref{improve}, 
we use methods of $L^2$-homology.
For a review of $L^2$-homology, see \cite{Lueck (1997)}.
Let $G$ and $\Gamma$ be as in the hypotheses of Theorem \ref{mainprop}.
Let $b_i^{(2)}(\Gamma) \in \R$ denote the $i$-th 
$L^2$-Betti number of $\Gamma$.
Let $\Rad$ be the radical of $G$, let $L$ be a Levi subgroup of $G$ and
let $K$ be the maximal compact connected
normal subgroup of $L$. Put $G_1 = \Rad \cdot K$ and
$G_2 = G/G_1$, a connected semisimple Lie group 
whose Lie algebra has no compact
factors. Let $\beta : G \rightarrow
G_2$ be the projection map. Put 
$\Gamma_1 = \Gamma \cap G_1$ and $\Gamma_2 = \beta(\Gamma)$.
Then there is an exact sequence
\begin{equation}
1 \longrightarrow \Gamma_1 \longrightarrow \Gamma 
\stackrel{\beta}{\longrightarrow} \Gamma_2 
\longrightarrow 1
\end{equation}
where $\Gamma_1$ is a lattice in $G_1$ and $\Gamma_2$ is a lattice
in $G_2$ \cite{Auslander (1963)}.

\begin{lemma} \label{b1}
If $b_1^{(2)}(\Gamma) \ne 0$ then  
$\Gamma$ has a finite normal subgroup $F$ such that $\Gamma/F$ is a 
lattice in $\PSL_2( \R)$.
\end{lemma}
\begin{pf}
There are the following possibilities:\\
{\bf A.} $\Gamma_1$ is infinite. Then $\Gamma$ has an infinite normal amenable
subgroup. By a result of Cheeger and Gromov, the $L^2$-Betti numbers of
$\Gamma$ vanish \cite[Theorem 10.12]{Lueck (1997)}.\\
{\bf B.} 
$\Gamma_1$ is finite and $\Gamma_2$ is finite. (That is, $\Gamma_2 = \{e\}$.) 
Then $\Gamma$ is finite
and $b_1^{(2)}(\Gamma) = 0$.\\
{\bf C.} 
$\Gamma_1$ is finite and $\Gamma_2$ is infinite. By the Leray-Serre spectral
sequence for $L^2$-homology,
$b_1^{(2)}(\Gamma) = b_1^{(2)}(\Gamma_2)/|\Gamma_1|$. Suppose that
$b_1^{(2)}(\Gamma_2) \ne 0$. If $G_2$ had an infinite center then 
$\Gamma_2$, being a lattice, would have to have an infinite center.  This
would imply by \cite[Theorem 10.12]{Lueck (1997)} that $b_1^{(2)}(\Gamma_2)$
vanishes, so $G_2$ must have a finite center $Z(G_2)$. Put $G_3 = G_2/Z(G_2)$,
let $\gamma : G_2 \rightarrow G_3$ be the projection and put 
$\Gamma_3 = \gamma(\Gamma_2)$, a lattice in $G_3$. 
Then there is the exact sequence
\begin{equation}
1 \longrightarrow \Gamma_2 \cap Z(G_2) \longrightarrow \Gamma_2
\stackrel{\gamma}{\longrightarrow} \Gamma_3 \longrightarrow 1
\end{equation}
and so $b_1^{(2)}(\Gamma_2) = b_1^{(2)}(\Gamma_3)/|\Gamma_2 \cap Z(G_2)|$.
Let $K_3$ be a maximal compact subgroup
of $G_3$ and let ${\cal F}$ be a fundamental domain for the $\Gamma_3$-action
on $G_3/K_3$. Let $\Pi(x,y)$ be the Schwartz kernel for the projection
operator onto the $L^2$-harmonic $1$-forms on $G_3/K_3$.
By \cite{Cheeger-Gromov (1985)}, 
$b_1^{(2)}(\Gamma_3) = \int_{\cal F} \tr(\Pi(x,x)) d\vol(x)$. Hence $G_3/K_3$
has nonzero $L^2$-harmonic $1$-forms.  By the K\"unneth formula for
$L^2$-cohomology and \cite[Section II.5]{Borel-Wallach (1980)},
the only possibility is $G_3 = \PSL_2( \R)$. Then there is the
exact sequence
\begin{equation}
1 \longrightarrow \Gamma \cap \Ker(\gamma \circ \beta) \longrightarrow \Gamma
\stackrel{\gamma \circ \beta}{\longrightarrow} \Gamma_3 \longrightarrow 1
\end{equation}
with $\Gamma \cap \Ker(\gamma \circ \beta)$ finite.
\end{pf}

Let $\gd \Gamma$ be the minimal dimension of a $K(\Gamma, 1)$-complex
\cite[p. 185]{Brown (1982)}. We will need the following result of Hillman
\cite[Theorem 2]{Hillman (1997)}.  For completeness, we give the short proof.

\begin{lemma} (Hillman) \label{Hillman}
If $\Gamma$ is a finitely-presented group then $\Def(\Gamma) \le
1 + b_1^{(2)}(\Gamma)$. Equality implies that there is a finite
$K(\Gamma, 1)$-complex $X$ with $\dim(X) \le 2$.
\end{lemma}
\begin{pf}
If $\Gamma$ is finite then $\Def(\Gamma) \le 0$, so we may assume that
$\Gamma$ is infinite.
Given a presentation of $\Gamma$ with $g$ generators and $r$ relations,
let $X$ be the corresponding $2$-complex. As $X$ is
two-dimensional, its second $L^2$-homology group is the same as the space of
square-integrable real cellular 
$2$-cycles on the universal cover
$\widetilde{X}$. This contains the ordinary integer cellular
$2$-cycles as a subgroup.

We have
\begin{equation}
\chi(X) = 1 - g + r = b_0^{(2)}(X) - b_1^{(2)}(X) + b_2^{(2)}(X) =
- b_1^{(2)}(\Gamma) + b_2^{(2)}(X).
\end{equation}
Hence 
\begin{equation} \label{ineq}
g - r = 1 + b_1^{(2)}(\Gamma) - b_2^{(2)}(X) \le 1 + b_1^{(2)}(\Gamma). 
\end{equation}
If $g - r = 1 + b_1^{(2)}(\Gamma)$ then $b_2^{(2)}(X) = 0$.
Hence $\HH_2(\widetilde{X}; \Z) = 0$. From the Hurewicz theorem, 
$\widetilde{X}$ is contractible.
\end{pf}

We now prove Theorem \ref{mainprop}. Suppose that
$\Def(\Gamma) > 0$. Then first of all, $|\Gamma| = \infty$.
Suppose that $\Gamma$ does not have
a finite normal subgroup $F$ such that $G/F$ is a lattice in $\PSL_2( \R)$.
By Lemma \ref{b1}, $b_1^{(2)}(\Gamma) = 0$. Then Lemma
\ref{Hillman} implies that $\Def(\Gamma) = 1$ and $\gd \Gamma \le 2$.
In particular, $\Gamma$ is torsion-free.
 
As $\Gamma_1$ is a lattice in $K \cdot \Rad$, it is a uniform lattice
\cite[Chapter III]{Raghunathan (1972)}.
Furthermore, as $\Gamma_1$ is a subgroup of $\Gamma$, $\gd \Gamma_1 \le 2$
and so $\Gamma_1$ must be $\{e\}$, $\Z$, $\Z^2$ or
the fundamental group of a Klein bottle. We go through the possibilities :\\
{\bf A.} $\Gamma_1 = \{e\}$. Then $\Gamma = \Gamma_2$ 
is a torsion-free lattice in the semisimple group $G_2$. 
Using a result of Borel and Serre \cite[p. 218]{Brown (1982)}, 
the fact that $\gd \Gamma \le 2$ implies that the Lie algebra of $G_2$ is
$\Sl_2( \R)$, $\Sl_2( \R) \oplus \Sl_2( \R)$ or $\Sl_2( \C)$.
One possibility is $G_2 = \widetilde{\PSL_2(\R)}$. Using the
embedding $\widetilde{\PSL_2(\R)} \cong \Z \times_{\Z} 
\widetilde{\PSL_2(\R)} \rightarrow \R \times_{\Z} 
\widetilde{\PSL_2(\R)}$, in this case we can say that $\Gamma$ is isomorphic 
to a lattice in $\R \times_{\Z} \widetilde{\PSL_2(\R)}$. On the other hand,
if $G_2$ is a finite covering of $\PSL_2(\R)$ then
$b_1^{(2)}(\Gamma) \ne 0$, contrary to assumption.
If $G_2$ is an infinite covering of $\PSL_2( \R) \times \PSL_2( \R)$ then the
Leray-Serre spectral sequence implies that
$\Gamma_2$ has cohomological dimension greater than two, contrary to
assumption. If $G_2$ is a finite covering of $\PSL_2( \R) \times \PSL_2( \R)$ 
then Lemma \ref{lattice} below will
show that $\Def(\Gamma) \le 0$, contrary to assumption. 
If $G_2 = \SL_2(\C)$, let $p : \SL_2(\C) \rightarrow \PSL_2(\C)$ be the
projection map.  Then there is the exact sequence
\begin{equation}
1 \longrightarrow \Gamma \cap \Ker(p) \longrightarrow \Gamma
\stackrel{p}{\longrightarrow} p(\Gamma) \longrightarrow 1.
\end{equation}
As $\Gamma$ is torsion-free, $\Gamma \cap \Ker(p) = \{e\}$ and so
$\Gamma$ is isomorphic to $p(\Gamma)$, a lattice in $\PSL_2(\C)$.  
Thus in any case, $\Gamma$ is
isomorphic to a torsion-free lattice in $\R \times_{\Z} 
\widetilde{\PSL_2( \R)}$ or
$\PSL_2( \C)$.
If $\Gamma$ is uniform then
$\gd \Gamma = 3$. Thus $\Gamma$ must be nonuniform. The torsion-free 
nonuniform lattices in $\R \times_{\Z} 
\widetilde{\PSL_2( \R)}$ and $\R \times
\PSL_2( \R)$ are isomorphic, as they both correspond to the
Seifert fiber spaces whose base is a hyperbolic orbifold with boundary 
\cite{Scott (1983)}. We conclude that $\Gamma$ is
isomorphic to a torsion-free nonuniform lattice in $\R \times \PSL_2( \R)$ or
$\PSL_2( \C)$.\\
{\bf B.} $\Gamma_1 = \Z$. Let $\Gamma_2^\prime$ be a finite-index torsion-free
subgroup of $\Gamma_2$ which acts trivially on $\Z$ 
and put $\Gamma^\prime = \beta^{-1}(\Gamma_2^\prime)$,
a finite-index subgroup of $\Gamma$. Then there is the exact sequence
\begin{equation}
1 \longrightarrow \Gamma_1 \longrightarrow \Gamma^\prime
\stackrel{\beta}{\longrightarrow} \Gamma_2^\prime \longrightarrow 1.
\end{equation}
Let $M$ be a $\Gamma_2^\prime$-module and let
$\beta^* M$ be the corresponding $\Gamma^\prime$-module.
If $\HH^*(\Gamma_2^\prime; M) \ne 0$,
let $k$ be the largest integer such that
$\HH^k(\Gamma_2^\prime; M) \ne 0$.
Then by the Leray-Serre spectral sequence,
$\HH^{k+1}(\Gamma^\prime; \beta^* M) \ne 0$. As
$\gd \Gamma^\prime \le 2$, we must have $k \le 1$. 
Thus the cohomological dimension of $\Gamma_2^\prime$ is at most one
and $\Gamma_2^\prime$ must be trivial or
a free group \cite[p. 185]{Brown (1982)}. 
If $\Gamma_2^\prime = \{e\}$ then $G_2 = \{e\}$ and $\Gamma = \Z$. If
$\Gamma_2^\prime$ is a free group then $G_2$ is a finite covering of
$\PSL_2( \R)$. Let $\sigma : G_2 \rightarrow \PSL_2( \R)$ be the
projection map and put $L = (\sigma \circ \beta)(\Gamma)$.  Then there is
the exact sequence
\begin{equation}
1 \longrightarrow \Gamma \cap \Ker(\sigma \circ \beta) \longrightarrow 
\Gamma
\stackrel{\sigma \circ \beta}{\longrightarrow} 
L \longrightarrow 1
\end{equation}
where $L$ is a lattice in $\PSL_2(\R)$ and 
$\Gamma \cap \Ker(\sigma \circ \beta)$ is
virtually cyclic.  As $\Gamma \cap \Ker(\sigma \circ \beta)$ is torsion-free, 
it must equal $\Z$. It follows that
$\Gamma$ is isomorphic to a lattice in 
$\R \times \PSL_2(\R)$ or $\R \times_{\Z} \widetilde{\PSL_2(\R)}$. 
If $\Gamma$ is uniform then
$\gd \Gamma = 3$. Thus $\Gamma$ is nonuniform and is isomorphic to a lattice in
$\R \times \PSL_2(\R)$. \\
{\bf C.} $\Gamma_1 = \Z^2$. Let $\Gamma_2^\prime$ be a finite-index 
torsion-free
subgroup of $\Gamma_2$ which acts on $\Z^2$ with determinant $1$ 
and put $\Gamma^\prime = \beta^{-1}(\Gamma_2^\prime)$,
a finite-index subgroup of $\Gamma$.
Let $M$ be a $\Gamma_2^\prime$-module and let
$\beta^* M$ be the corresponding $\Gamma^\prime$-module.
If $\HH^*(\Gamma_2^\prime; M) \ne 0$,
let $k$ be the largest integer such that
$\HH^k(\Gamma_2^\prime; M) \ne 0$.
Then by the Leray-Serre spectral sequence,
$\HH^{k+2}(\Gamma^\prime ; \beta^* M) \ne 0$. As
$\gd \Gamma^\prime \le 2$, we must have $k = 0$. Thus the cohomological
dimension of $\Gamma_2^\prime$ is zero, so $\Gamma_2^\prime = \{e\}$ and 
$G_2 = \{e\}$. Then $\Gamma = \Z^2$.\\
{\bf D.} $\Gamma_1$ is the fundamental group of a Klein bottle. Let $\Z^2$
be the unique maximal abelian subgroup of $\Gamma_1$. Any automorphism of
$\Gamma_1$ acts as an automorphism of $\Z^2$. Thus we get a homomorphism 
$\phi : \Aut(\Gamma_1) \rightarrow \GL_2(\Z)$. Let $\rho : \Gamma
\rightarrow \Aut(\Gamma_1)$ be given by $(\rho(\gamma))(\gamma_1) =
\gamma \gamma_1 \gamma^{-1}$.
Put $\widetilde{\Gamma} =
\Ker(\det \circ \phi \circ \rho)$, an index-$2$ subgroup of $\Gamma$, 
and put $\widetilde{\Gamma_2} = \beta \left(
\widetilde{\Gamma} \right)$. Then there is an exact sequence
\begin{equation}
1 \longrightarrow \Z^2 \longrightarrow \widetilde{\Gamma}
\stackrel{\beta}{\longrightarrow} \widetilde{\Gamma_2} \longrightarrow 1.
\end{equation}
As in case C, it follows that $G_2 = \{e\}$ and $\Gamma = \Gamma_1$ is the
fundamental group of a Klein bottle.

This proves Theorem \ref{mainprop}. We now prove Theorem \ref{improve}.
Let $X$ be as in the proof of Lemma \ref{Hillman}. As the classifying
map $X \rightarrow B\Gamma$ is $2$-connected, $b_2^{(2)}(X) \ge
b_2^{(2)}(\Gamma)$. Then from (\ref{ineq}),
\begin{equation} \label{defineq}
\Def(\Gamma) \le 1 + b_1^{(2)}(\Gamma) - b_2^{(2)}(\Gamma).
\end{equation}
For the lattices in question, let $G$ be the Lie group, let $K$ now be a
maximal compact subgroup of $G$ 
and put $M = \Gamma \backslash G/K$, an orbifold. 
As $G/K$ has no $L^2$-harmonic $1$-forms
\cite[Section II.5]{Borel-Wallach (1980)}, it follows from
\cite{Cheeger-Gromov (1985)} that $b_1^{(2)}(\Gamma) = b_3^{(2)}(\Gamma) = 0$.
As $|\Gamma| = \infty$, we have $b_0^{(2)}(\Gamma) = b_4^{(2)}(\Gamma) = 0$.
If $\chi(\Gamma)$
is the rational-valued group Euler characteristic of $\Gamma$ 
\cite[p. 249]{Brown (1982)} then
\begin{equation} \label{Euler}
\chi(\Gamma) =  b_0^{(2)}(\Gamma) -  b_1^{(2)}(\Gamma) +
b_2^{(2)}(\Gamma) - b_3^{(2)}(\Gamma) +
b_4^{(2)}(\Gamma) = b_2^{(2)}(\Gamma).
\end{equation}
From (\ref{defineq}) and (\ref{Euler}), we obtain
\begin{equation} \label{defineq2}
\Def(\Gamma) \le 1 - \chi(\Gamma).
\end{equation}
Furthermore, letting
$e(M,g) \in \Omega^4(M)$ denote the Euler density, it 
follows from \cite{Cheeger-Gromov (1985)} that 
\begin{equation} \label{integral}
\chi(\Gamma) = \int_M e(M,g).
\end{equation}
Let $G^d/K$ be the compact dual symmetric space to $G/K$. By the
Hirzebruch proportionality principle,
\begin{equation}
\frac{\int_M e(M,g)}{\chi(G^d/K)} =
\frac{\vol(M)}{\vol(G^d/K)}.
\end{equation}
We have the table
\begin{equation*}
\begin{array}{ccccccc}
\underline{G} & \Big| & \underline{G^d/K} & \Big| & 
\underline{\chi(G^d/K)} & \Big| & \underline{\vol(G^d/K)} \\
% & \Big| &  & \Big| &  & \Big| & \\
\SO(4,1) & \Big| & S^4 & \Big| & 2 & \Big| & \frac{8 \pi^2}{3} \\
\SU(2,1) & \Big| & \C P^2 & \Big| & 3 & \Big| & \frac{\pi^2}{2} \\
\PSL_2(\R) \times \PSL_2(\R) & \Big| & S^2 \times S^2 & \Big| & 4 & \Big| & 
16 \pi^2.
\end{array}
\end{equation*}
This proves Theorem \ref{improve}.

\begin{lemma} \label{lattice}
Let $G$ be a connected Lie group with a surjective homomorphism
$\rho : G \rightarrow \PSL_2(\R) \times \PSL_2(\R)$
such that $\Ker(\rho)$ is central in $G$ and finite.
If $\Gamma$ is a lattice in $G$ then $\Def(\Gamma) \le 0$.
\end{lemma}
\begin{pf}
Equation (\ref{defineq2}) is still valid for $\Gamma$.
We have $\chi(\Gamma) = 
\chi(\rho(\Gamma))/|\Gamma \cap \Ker(\rho)|$. Applying 
(\ref{integral}) to $\rho(\Gamma)$, 
the proof of Theorem \ref{improve} gives $\chi(\rho(\Gamma)) > 0$.
Hence $\chi(\Gamma) > 0$ and $\Def(\Gamma) \le 0$. 
\end{pf}

\end{document}